\documentclass[12pt,leqno]{amsart}
\usepackage[dvips]{graphicx}
\usepackage{amsfonts}
\usepackage{amssymb}
\usepackage{amsmath}
\usepackage{amscd}
\usepackage{graphicx}

\def\DATE{\today}
\newtheorem{theorem}{Theorem}
\newtheorem{definition}[theorem]{Definition}

\newtheorem{proposition}[theorem]{Proposition}

\newcommand\g{\mathfrak{g}}

\newcommand\K{\mathbb{K}}
\newcommand\Z{\mathbb{Z}}

\newcommand\p{\mathcal{P}}

\email{elisabeth.remm@uha.fr}

\title{Poisson Superalgebras as nonassociative algebras}
\author{Elisabeth Remm}
\date{6 mai 2012}
\address{Universit\'{e} de Haute Alsace, LMIA, 4 rue des Fr\`{e}res Lumi\`{e}re, 68093 Mulhouse}

\begin{document}

\maketitle

\begin{abstract}
Poisson superalgebras
 are known as a $\mathbb{Z}_2$-graded vector space with two operations, 
an associative supercommutative multiplication and a super bracket tied up by the super Leibniz relation. 
We show that we can consider a single nonassociative multiplication containing all these datas and then consider Poisson superalgebra
 as non associative algebras. 
\end{abstract}

\section{Generalities on Poisson algebras}

\subsection{Poisson algebras}
Let $\K$ be a  field of characteristic $0$. A $\K$-Poisson algebra is a $\K$-vector space  $\p$ equipped with two bilinear products denoted by $x\cdot y$  and $\{x ,y \}$, 
having the following properties:
\begin{enumerate}
\item The couple $(\p,\cdot)$ is an associative commutative $\K$-algebra.

\item The couple $(\p, \{ , \})$ is  a $\K$-Lie algebra.

\item The products $\cdot$ and $\{, \}$ satisfy the Leibniz rule:
$$\{x\cdot y,z\}=x \cdot\{y,z\}+\{x,z\}\cdot y,$$
for any $x,y,z \in \p.$
\end{enumerate}
The product $\{,\}$ is usually called Poisson bracket and the Leibniz identity means that the Poisson bracket acts as a derivation of the associative product.

In \cite{markl-eli-poisson}, one proves that any Poisson structure on a $\K$-vector space is also given by a nonassociative product, denoted by $xy$ and satisfying the non associative identity
\begin{eqnarray}
\label{associator} 3A(x,y,z)=(xz)y+(yz) x-(y x)z-(z x) y.
\end{eqnarray}
where $A(x,y,z)$ is the associator $A(x,y,z)=(xy)z-x(yz)$. In fact, if $\p$ is a Poisson algebra given by the associative product $x\cdot y$ and the Poisson bracket $\{x,y\}$, 
then $xy$ is given by $$xy=\{x,y\}+x \cdot y.$$ Conversely, the Poisson bracket and the associative product of $\p$ are the skew-symmetric part and the symmetric part of 
the product $xy$. Thus it is equivalent to present a Poisson algebra classically or by this nonassociative product. In \cite{Mic-Eli-Poisson}, we have studied algebraic 
properties of the nonassociative algebra $\p$. 
In particular we have proved that this algebra is flexible, power-associative and admits a Pierce decomposition.

\medskip

If $\p$ is a Poisson algebra given by the nonassociative product (\ref{associator}), 
we denote by $\g_{\p}$ the Lie algebra on the same vector space $\p$ whose Lie bracket is $$\{x,y\}=\displaystyle\frac{xy-yx}{2}$$ and by $\mathcal{A}_{\p}$ 
the commutative associative algebra, on the same vector space, whose product is $$x \cdot y=\displaystyle\frac{xy+yx}{2}.$$
An important problem in mathematical physics and more precisely in Quantum Field theory is the deformation of Poisson algebras. 
The classical deformations of Poisson algebras consist of deformations of the Poisson brackets, that is, deformations of $\g_\mathcal{P}$ which let the
associatif multiplication of $\mathcal{A}_\mathcal{P}$ unchanged and satisfying the Leibniz identity \cite{Pich}. In \cite{remm-deformPoisson} 
this type of deformation has been generalized  by using a nonassociative multiplication defining the Poisson stucture. The deformations of this nonassociative
multiplication provides general Poisson deformations.


\section{Poisson superalgebra
 }
By a $\K$-super vector space, we mean a $\Z_2$-graded vector space $V=V_0 \oplus V_1$. The vectors of $V_0$ and $V_1$ are called homogeneous vectors of degree respectively equal to $0$ and $1$. For an homogeneous vector $x$, we denote by $\mid x\mid$ its degree.
A  $\K$-Poisson superalgebra
 is a $\K$-super vector space  $\p=\p_0 \oplus \p_1$ equipped with two bilinear products denoted by $x\cdot y$  and $\{x ,y \}$, having the following properties:
\begin{itemize}
\item The couple $(\p,\cdot)$ is a associative super commutative $\K$-algebra,
that is, $$x \cdot y=(-1)^{\mid x \mid \mid y \mid}y \cdot x.$$ 
\item The couple $(\p, \{ , \})$ is  a $\K$-Lie super algebra, that is,
$$\{x ,y \}=-(-1)^{\mid x \mid \mid y \mid}\{y ,x \}$$
and satisfying the super Jacobi condition:
$$(-1)^{\mid z \mid \mid x \mid}\{x , \{y ,z \}\}+(-1)^{\mid x \mid \mid y \mid}\{y , \{z ,x \}\}+
(-1)^{\mid y \mid \mid z \mid}\{z , \{x ,y \}\}=0.$$
\item The products $\cdot$ and $\{, \}$ satisfy the super Leibniz rule:
$$\{x,y\cdot z\}=\{x,y\} \cdot z+(-1)^{\mid x \mid \mid y \mid}y \cdot \{x,z\}.$$
where $x,y$ and $z$ are homogeneous vectors.
\end{itemize}

\bigskip

\begin{theorem}
Let $\p$ a $\K$-super vector space. Thus $\p$ is a Poisson superalgebra
 if and only if there exists on $\p$ a nonassociative product $x  y$ satisfying
 \begin{equation}\label{superPoisson}
\left\{
 \begin{array}{l}
\medskip
 3(xy)z-3x(yz)  + (-1)^{\mid x \mid \mid y \mid}(yx)z -(-1)^{\mid y \mid \mid z \mid} (xz)y
 - (-1)^{\mid x \mid \mid y \mid+\mid x \mid \mid z \mid}(yz)x \\
 + (-1)^{\mid x \mid \mid z \mid+\mid y \mid \mid z \mid}(zx)y =0 
\end{array}
\right.
\end{equation}
for any homogeneous vectors $x,y,z \in\p$.
\end{theorem}

\noindent{\it Proof.} Assume that $(\p,\cdot,\{,\})$ is a Poisson superalgebra
.
Consider the multiplication
$$\begin{array}{l}
xy=x \cdot y + \{ x,y \}.
\end{array}$$
We deduce that
$$\displaystyle x \cdot y=\frac{1}{2}(xy+(-1)^{\mid x \mid \mid y \mid}yx).$$
Thus
the associativity condition writes for homogeneous vectors
$$\begin{array}{rl}
\medskip
v_1(x,y,z)= & A(x,y,z)-(-1)^{\mid x \mid \mid y \mid+\mid x \mid \mid z \mid+\mid y \mid \mid z \mid}A(z,y,x)
 +(-1)^{\mid x \mid \mid y \mid}(yx)z\\
 & -(-1)^{\mid y \mid \mid z \mid}x(zy)-(-1)^{\mid x \mid \mid y \mid+\mid x \mid \mid z\mid}(yz)x+(-1)^{\mid x \mid \mid z \mid+\mid y \mid \mid z \mid}z(xy)\\
 = &0
\end{array}$$
where $A(x,y,z)=(xy)z-x(yz)$. 
Likewise, the Poisson bracket writes for homogeneous vectors
$$\displaystyle \{x , y\}=\frac{1}{2}(xy-(-1)^{\mid x \mid \mid y \mid}yx)$$
and the super Jacobi condition
$$\begin{array}{rl}
\medskip
v_2(x,y,z)=& (-1)^{\mid x \mid \mid z \mid}A(x,y,z)-(-1)^{\mid x \mid \mid y \mid+\mid x \mid \mid z \mid}A(y,x,z)-(-1)^{\mid x \mid \mid y \mid+\mid y \mid \mid z \mid}A(z,y,x)\\
& -(-1)^{\mid x \mid \mid z \mid+\mid y \mid \mid z \mid}A(x,z,y)+(-1)^{\mid x \mid \mid y \mid}A(y,z,x)+(-1)^{\mid y \mid \mid z \mid}A(z,x,y) \\
=&0
\end{array}$$
The super Leibniz writes
$$\begin{array}{rl}
\medskip
v_3(x,y,z)=& A(x,y,z)-(-1)^{\mid x \mid \mid y \mid}A(y,x,z)+(-1)^{\mid x \mid \mid y \mid+\mid x \mid \mid z \mid+\mid y \mid \mid z \mid}A(z,y,x)\\
& +(-1)^{\mid y \mid \mid z \mid}A(x,z,y)
 +(-1)^{\mid x \mid \mid y \mid+\mid x \mid \mid z \mid}A(y,z,x)-(-1)^{\mid x \mid \mid z \mid+\mid y \mid \mid z \mid}A(z,x,y) \\
 =& 0.
\end{array}$$
Let us consider the vector

$$\begin{array}{ll}
\medskip
v(x,y,z)  = 
 & \frac{1}{3}\left[ (-1)^{\mid x \mid \mid y \mid}(yx)z -(-1)^{\mid y \mid \mid z \mid} (xz)y
 - (-1)^{\mid x \mid (\mid y \mid +\mid z \mid}(yz)x
  + (-1)^{(\mid x \mid +\mid y \mid )\mid z \mid}(zx)y\right] \\
& +(xy)z-x(yz) .\\
\end{array}$$
Then 
$$\begin{array}{l}
v(x,y,z)= \frac{1}{6}\left( 2v_1(x,y,z)+(-1)^{\mid x \mid \mid z \mid}v_2(x,y,z)+v_3(x,y,z)+2(-1)^{\mid x \mid \mid z \mid+\mid y \mid \mid z \mid}v_3(z,x,y) \right).
\end{array}
$$
We deduce that the product $xy$ satisfies
$$v(x,y,z)=0$$
for any homogeneous vectors $x,y,z$.

\noindent Conversely, assume that the product of the non associative product $\p$ satisfies $v(x,y,z)=0$ for any homogeneous vestors $x,y,z.$ Let 
$v_1(x,y,z), v_2(x,y,z), v_3(x,y,z)$ be the vectors of $\p$ defined in the first part respectively in relation with the associativity, the super Jacobi and super Leibniz relations. 
We have
$$\begin{array}{ll}
\medskip
v_1(x,y,z)=&v(x,y,z)-(-1)^{\mid x \mid \mid y \mid+\mid x \mid \mid z \mid+\mid y \mid \mid z \mid}v(z,y,x)
+(-1)^{\mid y \mid \mid z \mid}v(x,z,y)\\
\medskip
& -(-1)^{\mid x \mid \mid z \mid+\mid y \mid \mid z \mid}v(z,x,y)\\
\medskip
v_2(x,y,z) =&(-1)^{\mid x \mid \mid z \mid}v(x,y,z)-(-1)^{\mid x \mid \mid y \mid+\mid x \mid \mid z \mid}v(y,x,z)-(-1)^{\mid x \mid \mid y \mid+\mid y \mid \mid z \mid}v(z,y,x)\\
\medskip
& -(-1)^{\mid x \mid \mid z \mid+\mid y \mid \mid z \mid}v(x,z,y)+(-1)^{\mid x \mid \mid y \mid}v(y,z,x)+(-1)^{\mid y\mid \mid z \mid}v(z,x,y)\\
\medskip
v_3(x,y,z)=&v(x,y,z)-(-1)^{\mid x \mid \mid y \mid}v(y,x,z)+(-1)^{\mid x \mid \mid y \mid+\mid x \mid \mid z \mid+\mid y \mid \mid z \mid}v(z,y,x)\\
& +(-1)^{\mid y \mid \mid z \mid}v(x,z,y)+(-1)^{\mid x \mid \mid y \mid+\mid x \mid \mid z \mid}v(y,z,x)-(-1)^{\mid x\mid \mid z \mid+\mid y\mid \mid z \mid}v(z,x,y)
\end{array}$$

\medskip
\noindent{\bf Examples.} Any $2$-dimensional superalgebra $\p=V_0 \oplus V_1$ with an homogeneous basis $\{e_0,e_1\}$ is defined
$$\left\{
\begin{array}{l}
e_0e_0=ae_0,\\
 e_0e_1=be_1, \ e_1e_0=ce_1,\\
e_1e_1=de_0.
\end{array}
\right.
$$
This is a super Poisson  multiplication if and only if we have
$$\left\{
\begin{array}{l}
d=0,\\
3(a-b)b+ab-2bc+c^2=0,\\
3(a-c)c+ab-2bc+c^2=0,
\end{array}
\right.
$$
or
$$\left\{
\begin{array}{l}
a=0,\\
a=b=c.
\end{array}
\right.
$$
We obtain the following $2$-dimensional Poisson superalgebras
$$\left\{
\begin{array}{ll}
\mathcal{SP}_{2,1} & e_0e_0=ae_1\\
\mathcal{SP}_{2,2} & e_0e_0=ae_1, e_0e_1=e_1e_0=ae_1\\
\mathcal{SP}_{2,3} &  e_0e_1=-e_1e_0=be_1\\
\mathcal{SP}_{2,4} &  e_1e_1=de_0,\\
\end{array}
\right.
$$
the non written product being considered equal to $0$.
Let us note that these $2$-dimensional algebras correspond to the algebras $(\mu_{16},\beta_2=0),(\mu_{16},\beta_2=1),(\mu_{9},\alpha_2=0,\beta_4=0),(\mu_{5},\alpha_2=0)$ in the classification \cite{goze-remm-2algebras}.
\section{Properties of Poisson superalgebras }
\begin{definition}
A nonassociative superalgebra  is called super flexive if  the multiplication $xy$ satisfy
$$A(x,y,z) + (-1)^{(|x||z|+|x||y|+|y||z|)}A(z,y,x)=0$$
for any homogeneous elemnts $x,y,z$, where $A(x,y,z)=(xy)z-x(yz)$ is the associator of the multiplication.
\end{definition}
\begin{proposition}
Let $\p$ be a Poisson superalgebra. Then the non associative product defining the super Poisson structure is super flexive.
\end{proposition}
\noindent{\it Proof.} In fact, let 
$$B(x,y,z)=3(A(x,y,z) + (-1)^{(|x||z|+|x||y|+|y||z|)}A(z,y,x)) .
$$
We  have
$$
\begin{array}
{rl}
\medskip
B(x,y,z) =& -(-1)^{\mid x \mid \mid y \mid}(yx)z +(-1)^{\mid y \mid \mid z \mid} (xz)y
 + (-1)^{\mid x \mid \mid y \mid+\mid x \mid \mid z \mid}(yz)x 
 - (-1)^{\mid x \mid \mid z \mid+\mid y \mid \mid z \mid}(zx)y \\
 \medskip
 
 &+ (-1)^{(|x||z|+|x||y|+|y||z|)}(-(-1)^{\mid z \mid \mid y \mid}(yz)x +(-1)^{\mid y \mid \mid x \mid} (zx)y \\
 \medskip
 
 &+ (-1)^{\mid z \mid \mid y \mid+\mid z \mid \mid x \mid}(yx)z 
 - (-1)^{\mid z \mid \mid x \mid+\mid y \mid \mid x \mid}(xz)y )\\
 \medskip
 
 =& (-(-1)^{\mid x \mid \mid y \mid} +(-1)^{\mid x \mid \mid y \mid})(yx)z + ((-1)^{\mid y \mid \mid z \mid} -(-1)^{\mid y \mid \mid z \mid}) (xz)y \\
 \medskip
 &+((-1)^{\mid x \mid \mid y \mid+\mid x \mid \mid z \mid} - (-1)^{\mid x \mid \mid y \mid+\mid x \mid \mid z \mid})(yz)x 
 +( - (-1)^{\mid x \mid \mid z \mid+\mid y \mid \mid z \mid}\\
 \medskip
 & +(-1)^{\mid x \mid \mid z \mid+\mid y \mid \mid z \mid})(zx)y \\
 \medskip
 =& 0.
 \end{array}
 $$
 
 \medskip
 
 \noindent{\bf Remark : On the power associativity.} Recall that a nonassociative algebra is power associative if every element
generates an associative subalgebra. Let $\p$ be a Poisson superalgebra
 provided with its non associative product $xy$. If $V_0$ is its the even homogeneous part, then the 
restriction of the product $xy$ is a multiplication in 
this homogeneous vector space satisfying Identity (\ref{superPoisson}). Since all the vectors of $V_0$ are of degree $0$, Identity (\ref{superPoisson}) is reduced to Identity
 (\ref{associator}). We deduce that $V_0$ is a Poisson algebra and any vector $x$ in $V_0$ generates an associative subalgebra of $V_0$ and of $\p$.
  
  Assume now that $y$ is an odd vector. We have
  $$y\cdot y=\displaystyle\frac{1}{2}(yy+(-1)yy)=0,$$
  and
  $$\{y,y\}=\displaystyle\frac{1}{2}(yy-(-1)yy)=yy.$$
  If we write $y^2=yy$, then $$y^2=\{y,y\}.$$
  This implies
  $$yy^2=y\{y,y\}=y \cdot \{y,y\}+\{y,\{y,y\}\}.$$
  But from the super identity of Jacobi, $\{y,\{y,y\}\}=0.$ Thus we have
  $$yy^2= y \cdot \{y,y\}=\{y,y\} \cdot y=y^2y.$$
  We can write
  $$y^3=yy^2=y^2y.$$
  Now
  $$y^2y^2=\{y,y\}\{y,y\}=\{y,y\}\cdot\{y,y\}+\{\{y,y\},\{y,y\}\}.$$
  We have also
  $$yy^3=y\cdot y^3+\{y,y^3\}=y\cdot y\cdot\{y,y\}+\{y,y\cdot \{y,y\}\}.$$
  But $y \cdot y=0$. Thus, from the Leibniz rule,
  $$yy^3=\{y,y\cdot \{y,y\}\}=-y\cdot \{y,\{y,y\}\}+\{y,y\}\cdot\{y,y\}=\{y,y\}\cdot \{y,y\}.$$
  We deduce
  $$y^2y^2-yy^3= \{\{y,y\},\{y,y\}\}.$$
  Since $\{y,y\}$ is of degree $0$, we obtain
$$y^2y^2-yy^3=0.$$
We can write
$$y^4=y^2y^2=yy^3=y^3y$$
the last equality results of $\{y,y\cdot \{y,y\}\}=\{y\cdot \{y,y\},y\}.$ Now, using  Identity (\ref{superPoisson}) to the triple $(y^i, y^j,y^k)$ with $i+j+k=5$, we obtain a 
linear system 
on the vectors $y^iy^j$ with $i+j=5$, which admits as solutions
$$yy^4=y^2y^3= y^3y^2=y⁴y.$$
Thus $y^5$ is well determinated. By induction, using Identity (\ref{superPoisson}) on the triple $(y^i, y^j,y^k)$ with $i+j+k=n$, using induction hypothesis $y^py^{n-1-p}=y^{n-1}$, we obtain that
$$y^n=y^py^{n-p}$$
for any $p=1,\cdots,n-1$. Thus any homogeneous element of edd degree generates an associative algebra.

\end{document}